\documentclass[12pt]{article}
\usepackage{mathrsfs}
\usepackage{amsmath,amsfonts,amssymb,rotating,amsthm}
\usepackage{hyperref}

\def\qed{\nopagebreak\hfill{\rule{4pt}{7pt}}}
\def\proof{\noindent {\it{Proof.} \hskip 2pt}}

\newtheorem{theo}{Theorem}[section]

\newtheorem{lemm}[theo]{Lemma}

\newtheorem{coro}[theo]{Corollary}
\newtheorem{conj}[theo]{Conjecture}

\theoremstyle{remark}

\newdimen\Squaresize \Squaresize=11pt
\newdimen\Thickness \Thickness=0.7pt
\def\Square#1{\hbox{\vrule width \Thickness
   \vbox to \Squaresize{\hrule height \Thickness\vss
    \hbox to \Squaresize{\hss#1\hss}
   \vss\hrule height\Thickness}
\unskip\vrule width \Thickness} \kern-\Thickness}

\def\Vsquare#1{\vbox{\Square{$#1$}}\kern-\Thickness}

\def\moins{\raise 1pt\hbox{{$\scriptstyle -$}}}

\begin{document}

\begin{center}
{\large \bf Zeta Functions and the Log-behavior of

 Combinatorial Sequences}
\end{center}

\begin{center}
William Y. C. Chen$^{1}$, Jeremy J. F. Guo$^{2}$ and Larry X. W. Wang$^{3}$ \\[8pt]
$^{1,2}$Center for Applied Mathematics\\
Tianjin University\\
 Tianjin 300072, P. R. China\\[6pt]
$^{3}$Center for Combinatorics, LPMC-TJKLC\\
Nankai University\\
 Tianjin 300071, P. R. China\\[6pt]

$^{1}${\tt chen@nankai.edu.cn}, $^{2}${\tt
guojf@mail.nankai.edu.cn}, $^{3}${\tt wsw82@nankai.edu.cn}
\end{center}

\vspace{0.3cm} \noindent{\bf Abstract.} In this paper, we use the
Riemann zeta function $\zeta(x)$ and the Bessel zeta function
$\zeta_{\mu}(x)$ to study the log-behavior of combinatorial
sequences. We prove that  $\zeta(x)$ is log-convex for $x>1$. As a
 consequence, we deduce that the sequence $\{|B_{2n}|/(2n)!\}_{n\geq 1}$ is log-convex,
 where $B_n$ is the $n$-th Bernoulli number.
 We introduce the function $\theta(x)=(2\zeta(x)\Gamma(x+1))^{\frac{1}{x}}$,
 where $\Gamma(x)$ is the gamma function, and we show that
$\log \theta(x)$ is strictly increasing for $x\geq 6$.
 This confirms a conjecture of Sun stating that the sequence
  $\{\sqrt[n] {|B_{2n}}|\}_{n\geq 1}$ is strictly increasing.
Amdeberhan, Moll and Vignat defined the numbers
$a_n(\mu)=2^{2n+1}(n+1)!(\mu+1)_n\zeta_{\mu}(2n)$ and conjectured
that  the sequence  $\{a_n(\mu)\}_{n\geq 1}$ is log-convex for
$\mu=0$ and $\mu=1$. By proving that $\zeta_{\mu}(x)$ is log-convex
 for $x>1$ and $\mu>-1$, we show that the sequence
$\{a_n(\mu)\}_{n\geq 1}$ is log-convex for any $\mu>-1$.
We introduce another function $\theta_{\mu}(x)$ involving
$\zeta_{\mu}(x)$ and the gamma function $\Gamma(x)$ and we show that
$\log \theta_{\mu}(x)$ is strictly increasing for $x>8e(\mu+2)^2$.
This implies that $\sqrt[n]{a_n(\mu)}<\sqrt[n+1]{a_{n+1}(\mu)}$ for
$n> 4e(\mu+2)^2$. Based on Dobinski's formula, we prove that
$\sqrt[n]{B_n}<\sqrt[n+1]{B_{n+1}}$ for $n\geq 1$, where $B_n$ is
the $n$-th
 Bell number. This confirms another conjecture of Sun.
 We also establish a connection between
   the increasing property of $\{\sqrt[n]{B_n}\}_{n\geq 1}$  and
   H\"{o}lder's inequality in probability theory.


\section{Introduction}

The objective of this paper is to present an analytic approach to
the log-behavior of   combinatorial sequences.

Let $B_n$ denote the $n$-th
Bernoulli number,  see  \cite{rota} and \cite{smi}.
Recall that $B_{2n+1}=0$ for $n\geq 1$ and $B_{2n}$ alternate in
sign for $n\geq 1$. We consider the log-behavior of the sequence
$\{|B_{2n}|\}_{n\geq 1}$.
A sequence $\{a_n\}_{n\geq 1}$ of real numbers
is said to be log-convex if for  $n\geq 2$,
\[
a_n^2\leq a_{n-1}a_{n+1}.
\]
It is well-known that
\begin{equation}\label{bernoulli}
\zeta(2n)=\frac{2^{2n-1}\pi^{2n}}{(2n)!}|B_{2n}|,
\end{equation}
where  \[ \zeta(x)=\sum_{n=1}^{\infty}\frac{1}{n^x}\]
is the Riemann zeta function. By proving that $\zeta(x)$ is
log-convex for $x>1$, we establish the log-convexity of the
sequence $\{|B_{2n}|/(2n)!\}_{n\geq 1}$.
Consequently, the sequence $\{|B_{2n}|\}_{n\geq 1}$ is log-convex.
 Moreover, we introduce the function
\begin{equation}\label{th}
\theta(x)  =(2\zeta(x)\Gamma(x+1))^{\frac{1}{x}},
\end{equation}
where  $\Gamma(x)$ is the gamma function.
We show that $\log \theta(x)$ is strictly increasing for
$x\geq 6$. From relation \eqref{bernoulli}, it can be seen that
\[ \sqrt[n]{|B_{2n}|}=\frac{1}{4\pi^2}\theta^2(2n).\]
So we reach the assertion that the sequence $\{\sqrt[n]
{|B_{2n}|}\}_{n\geq 1}$
 is strictly increasing. This confirms  a conjecture of Sun \cite{sun}, which
has been independently proved by  Luca and St\u{a}nic\u{a} \cite{luca}.
We conjecture that $(\log \theta(x))^{\prime \prime}<0$ for $x\geq 6$.

Our approach also applies to the sequence of generalized Lasalle
numbers. Let $C_n$ denote the $n$th Catalan number, that is,
$$C_n=\frac{1}{n+1}{{2n}\choose n},$$
and let $N_r(z)$ denote the $r$-th Narayana polynomial as given by
$$N_r(z)=\sum \limits_{k=1}^{r}\frac{1}{r}{r\choose {k-1}}{r\choose k}z^k.$$
 Lasalle \cite{las} derived the recurrence relation
\[
(z+1)N_r(z)-N_{r+1}(z)=\sum \limits_{n\geq 1}(-z)^n{{r-1}
\choose {2n-1}}A_nN_{r-2n+1}(z),
\]
where the numbers $A_n$ satisfy the recurrence relation
\begin{equation}\label{A}
  (-1)^{n-1}A_n = C_n + \sum_{j=1}^{n-1}(-1)^{j}{{2n-1}\choose{2j-1}}A_jC_{n-j}.
\end{equation}
Let \[ a_n={ 2A_n\over C_n}.\] Lasalle \cite{las} showed that
$\{a_n\}_{n\geq 1}$ is an increasing sequence of  positive integers.
  Amdeberhan, Moll and Vignat \cite{amd} established a
connection between $a_n$ and the Bessel zeta functions
$\zeta_{\mu}(x)$. Recall that for a real number $\mu$, the Bessel
function $J_{\mu}(z)$ of the first kind of order $\mu$ is defined by
\[
 J_{\mu}(z)=\left(\frac{z}{2}\right)^{\mu}\sum_{k=0}^{\infty}\frac{(-1)^k}
 {\Gamma(\mu+k+1)k!}\left(\frac{z}{2}\right)^{2k}.
\]
For $\mu\geq -1$, $J_{\mu}(z)$ has infinitely many positive real
zeros $j_{\mu, n}$, where we assume that
\[
  0\ <\ j_{\mu,1}\ <\ j_{\mu,2}\ <\ j_{\mu,3}\ <\ \cdots,
\]
see \cite[Sect. 4.14]{and}. The Bessel zeta functions $\zeta_{\mu}(x)$
are defined by
\begin{equation}\label{besselzeta}
  \zeta_{\mu}(x)=\sum_{n=1}^{\infty}\frac{1}{j_{\mu, n }^x}.
\end{equation}
Amdeberhan, Moll and Vignat \cite{amd} found the following relation
\begin{equation}\label{bessel2}
  a_n=2^{2n+1}(n+1)!(n-1)!\zeta_1(2n).
\end{equation}
They also gave the following generalization of $a_n$ for $\mu\geq
-1$,
\begin{equation}\label{a}
  a_n(\mu)= 2^{2n+1}(n-1)!(\mu +1)_n\zeta_{\mu}(2n),
\end{equation}
where $(\mu+1)_n= (\mu+1)(\mu+2)\cdots(\mu+n)$.

It is easily seen that $a_n=a_n(1)$. Setting $\mu=0$ in \eqref{a},
Amdeberhan, Moll and Vignat defined the sequence $\{b_n\}_{n\geq 1}$
as given by
\begin{equation}\label{bessel3}
b_n=\frac{1}{2}a_n(0)=2^{2n}n!(n-1)!\zeta_0(2n).
\end{equation}
Note that this sequence has  been studied by  Carlitz \cite{car}. It
is listed as Sequence A002190 in \cite{slo}.

Amdeberhan, Moll and Vignat conjectured that the sequences
$\{a_n\}_{n\geq 1}$  and $\{b_n\}_{n\geq 1}$ are log-convex. We show
that $\zeta_{\mu}(x)$ is log-convex for $x>1$. This implies that the
sequence  $\{a_n(\mu)\}_{n\geq 1}$ is log-convex for any $\mu>-1$.
This confirms the above conjectures, which have been independently
proved by Wang and Zhu \cite{wang3}.

Moreover, we define the following function
\begin{equation}\label{thetamu}
  \theta_{\mu}(x)=\left(\frac{2}{\mu !}\Gamma\left(\frac{x}{2}\right)
  \Gamma\left(\frac{x}{2}+\mu+1\right)\zeta_{\mu}(x)\right)^{\frac{1}{x}}.
\end{equation}
It can be easily checked that \begin{equation}
4\theta_{\mu}^2(2n)=\sqrt[n]{a_n(\mu)}.
 \end{equation}
We show that $\log \theta_{\mu}(x)$ is strictly increasing for
$x>8e(\mu+2)^2$. This leads to the increasing property  that
\begin{equation}\label{lasalle1}
\sqrt[n]{a_n(\mu)}<\sqrt[n+1]{a_{n+1}(\mu)},
\end{equation}
 for $n> 4e(\mu+2)^2$. We note that for $\mu=0$ and $\mu=1$ the
above relation \eqref{lasalle1} has been independently proved by
Wang and Zhu \cite{wang3}.

Owing to the formula of Dobinski, we may use our analytic approach
to study the log-behavior of Bell numbers.
 Let $B_n$ be
the $n$-th Bell number, that is, the number of partitions of $\{1,\
2,\ \ldots,\ n\}$, see \cite{bel} and \cite{rota2}. Notice that
 we have adopted  the same notation $B_n$ for both Bell numbers and
 Bernoulli
 numbers.
Recall that Dobinski's formula for the Bell numbers states that
\begin{equation*}
  B_n=\frac{1}{e}\sum_{k=0}^{\infty}\frac{k^n}{k!}.
\end{equation*}
For $x>0$, we define
\begin{equation}\label{bell}
  B(x)=\frac{1}{e}\sum_{k=0}^{\infty}\frac{k^x}{k!},
\end{equation}
so that  we have $B_n=B(n)$ whenever $n$ is a nonnegative integer.

 We show that $\log B(x)^{1/x}$ is increasing for $x\geq 1$. This
 implies that the sequence $\{\sqrt[n]{B_n}\}_{n\geq 1}$ is increasing, as
 conjectured by Sun \cite{sun}. We  conjecture that
 $(\log B(x)^{1/x})^{\prime\prime}<0$ for $x\geq 1$. In the last
  section, we give a probabilistic proof of the increasing property
  of the sequence $\{\sqrt[n]{B_n}\}_{n\geq 1}$ by using H\"{o}lder's
  inequality.

\section{The log-convexity of Bernoulli numbers}

To prove the log-convexity of Bernoulli numbers, we consider the
log-behavior of the Riemann zeta function
$\zeta(x)$ for $x>1$. Recall that a positive function $f$ is called
log-convex on a real interval $I=[a,b]$, if for all $x,y\in [a,b]$
and $\lambda \in [0,1]$,
\begin{equation}\label{df}
f(\lambda x+(1-\lambda)y)\leq f(x)^{\lambda}f(y)^{1-\lambda},
\end{equation}
see, for example, Artin \cite{art}. It is known that  a positive
function $f$ is log-convex if and only if $(\log f(x))^{\prime
\prime}\geq 0$. So, if
\begin{equation}\label{lr}
(\log \zeta(x))^{\prime \prime}>0,
\end{equation} for $x>1$, then we can deduce that
  $\zeta(x)$ is log-convex for $x>1$.

\begin{lemm}\label{riemann}
The Riemann zeta function $\zeta(x)$ is log-convex for $x>1$.
\end{lemm}

\proof Clearly,  condition
 \eqref{lr} is equivalent to
\begin{equation}\label{ri}
\zeta(x)\cdot \zeta^{\prime\prime}(x)-(\zeta^{\prime}(x))^2>0.
\end{equation}
Since $\zeta(x)$ converges for $x>1$,
 we find that  for $x>1$,
\begin{eqnarray*}
 \lefteqn{{\hspace{-0.3cm}\zeta(x)\zeta^{\prime \prime}(x)-(\zeta'(x))^2}}\\[6pt]
&=&\sum_{m=1}^{\infty}\frac{1}{m^x}\sum_{n=1}^{\infty}\frac{(\log
n)^2}{n^x}
-\sum_{m=1}^{\infty}\frac{\log m}{m^x}\sum_{n=1}^{\infty}\frac{\log n}{n^x}\\[6pt]
&=&\sum_{n>m\geq 1}\frac{(\log n)^2+(\log m)^2-2\log m \log n}{(mn)^x}\\[6pt]
&=&\sum_{n>m\geq 1}\frac{(\log n-\log m)^2 }{(mn)^x},
\end{eqnarray*}
 which is positive.  This completes the proof.
\qed

The log-convexity of $\zeta(x)$ enables us to deduce the following
  property of Bernoulli numbers.

\begin{theo}\label{str}
The sequence $\left\{\frac{|B_{2n}|}{(2n)!}\right\}_{n\geq 1}$ is
log-convex.
\end{theo}

\proof
Since $\zeta(x)$ is log-convex, setting $x=2n-2$, $y=2n+2$ and
$\lambda=1/2$ in the defining relation \eqref{df}, we find that
\begin{equation}\label{rn}
\zeta(2n-2)\zeta(2n+2)\geq \zeta(2n)^2.
\end{equation}
Invoking relation \eqref{bernoulli} between $\zeta(x)$ and $B_n$, we
obtain that
\[
\left(\frac{|B_{2n}|}{(2n)!}\right)^2\leq \frac{|B_{2n-2}|}{(2n-2)!}\cdot
 \frac{|B_{2n+2}|}{(2n+2)!}.
\]
This completes the proof.
\qed

Since $((2n)!)^2<(2n-2)!\cdot (2n+2)!$ for  $n\geq 1$,
 the above
theorem implies the following property.

\begin{coro}
The sequence $\{|B_{2n}|\}_{n\geq 1}$ is log-convex.
\end{coro}

%
%

\section{The log-behavior of $\theta(x)$}


In this section, we consider the log-behavior of the function
\[
\theta(x)=(2\zeta(x)\Gamma(x+1))^{\frac{1}{x}}.
\]
We begin with the following monotone property of $\log \theta(x)$.

\begin{theo}\label{rg}
$\log \theta(x)$ is  strictly increasing for $x\geq 6$.
\end{theo}

\proof To prove that $\log \theta(x)$ is increasing for $x\geq 6$,
we
 aim to show that
\begin{equation}\label{m1}
(\log \theta(x))^{\prime}>0,
\end{equation}
for $x\geq 6$. Let \[ g(x)=2\zeta(x)\Gamma(x+1).\] Then we have \[
\theta(x)=g(x)^{1/x}\] and
\[
(\log \theta(x))^{\prime}=\frac{1}{x}\left(\frac{g^{\prime}(x)}{g(x)}
-\frac{\log g(x)}{x}\right).
\]
Thus  \eqref{m1} can be rewritten as
\[
\frac{g^{\prime}(x)}{g(x)}>\frac{\log g(x)}{x},
\]
for $x\geq 6$. Since $\zeta(x)$ and $\Gamma(x)$ are continuous and
differentiable on $(1, \infty)$, so is $g(x)$  on $(1,\infty)$.
Applying the mean value theorem to $\log g(x)/x$, it can be shown
that there exists $t$ in
 $(2, x)$ such that
\begin{equation}\label{m2}
\frac{g(t)^{\prime}}{g(t)}>\frac{\log g(x)}{x}.
\end{equation}
Since $\zeta(2)=\frac{\pi^2}{6}$ and $\Gamma(3)=2$, we find that
\begin{equation}\label{lng2}
\log g(2)=\log (2\zeta(2)\Gamma(3))=\log \frac{2\pi^2}{3}<2.
\end{equation}
On the other hand, for $x\geq 6$, it is easily seen that $\zeta(x)>1$
and $\Gamma(x+1)>e^x$. It follows that
\begin{equation}\label{lng6}
\log g(x)=\log 2+\log \zeta(x)+\log \Gamma(x+1)>x.
\end{equation}
In  view of  \eqref{lng2} and \eqref{lng6},  we deduce  that for
$x\geq 6$,
\begin{equation}\label{m4}
 \frac{\log g(x)}{x}=\frac{(1-2/x)\log g(x)}{(1-2/x)x}<
 \frac{\log g(x)-2}{x-2}<\frac{\log g(x)-\log g(2)}{x-2}.
\end{equation}
Applying the mean value theorem to $\log g(x)$,
 we see that there exists $t\in (2,x)$ such that
\begin{equation}\label{m5}
(\log g(t))^{\prime}=\frac{\log g(x)-\log g(2)}{x-2},
\end{equation}
that is,
\begin{equation}\label{m6}
\frac{g^{\prime}(t)}{g(t)}=\frac{\log g(x)-\log g(2)}{x-2}.
\end{equation}
Combining \eqref{m4} and \eqref{m6}, we get \eqref{m2}.

Now we proceed to show that
\begin{equation}\label{m3}
\frac{g(x)^{\prime}}{g(x)}>\frac{g(t)^{\prime}}{g(t)}.
\end{equation}
Clearly, \eqref{m3} is equivalent to
\begin{equation}\label{m7}
\left(\frac{g^{\prime}(y)}{g(y)}\right)^{\prime}>0.
\end{equation}
By the definition of $g(x)$, we have
\[
\left(\frac{g^{\prime}(y)}{g(y)}\right)^{\prime}=(\log g(y))^{\prime \prime}=
(\log \Gamma(y+1))^{\prime \prime}+(\log \zeta(y))^{\prime \prime}.
\]
It is known that $(\log \Gamma(y+1))^{\prime \prime}>0$ for $y>1$,
see Andrews, Askey and Roy \cite[Theorem. 1.2.5]{and}. On the other
hand, in the proof of Lemma \ref{riemann}, we have shown that $(\log
\zeta(y))^{\prime \prime}>0$. This proves (\ref{m7}). In other
words,  $\frac{g^{\prime}(y)}{g(y)}$ is strictly increasing for
$y>1$. Thus for $2<t<x$,  inequality \eqref{m3} holds.

Combining \eqref{m2} and \eqref{m3}, we deduce that for $x\geq 6$,
\[
 \frac{g^{\prime}(x)}{g(x)}-\frac{\log g(x)}{x}>
 \frac{g^{\prime}(x)}{g(x)}-\frac{g^{\prime}(t)}{g(t)}>0.
\]
Thus $(\log \theta(x))^{\prime}>0$ for $x\geq 6$. This completes the
proof. \qed

From the log-behavior of $\theta(x)$, we are led to
an affirmative answer to a conjecture of Sun \cite{sun}.

\begin{coro}
The sequence $\{\sqrt[n]{|B_{2n}|}\}_{n\geq 1}$ is strictly increasing.
\end{coro}

\proof From  relation \eqref{bernoulli}, we see that for $n\geq 1$,
\begin{equation}\label{c1}
\sqrt[n]{|B_{2n}|}=\frac{1}{4\pi^2}\sqrt[n]{2\zeta(2n)(2n)!}
=\frac{1}{4\pi^2}\theta^2(2n).
\end{equation}
Since $\log \theta(x)$ is strictly increasing for $x \geq 6$,
 we find that $\theta(x)$ is also strictly increasing for  $x\geq 6$.
 It follows from \eqref{c1} that $\sqrt[n]{|B_{2n}|}$ is strictly
 increasing for $n\geq 3$. On the other hand, it is easily checked that
\[
 |B_2|< \sqrt{|B_4|}< \sqrt[3]{|B_6|}.
\]
This completes the proof.\qed

 The  conjecture of Sun was independently
proved by Luca and St\u{a}nic\u{a} \cite{luca}. In fact, they proved
that the sequence $\{\sqrt[n]{|B_{2n}|}\}_{n\geq 1}$ is log-concave,
which was also conjectured by Sun \cite{sun}.

We pose the following conjecture concerning the function
$\theta(x)$. If it is true, then it implies that the sequence
$\{\sqrt[n]{|B_{2n}|}\}_{n\geq 1}$ is log-concave.

\begin{conj}
The function $\theta(x)$ is log-concave for $x\geq 6$, that is,
for $x\geq 6$, $(\log f(x))^{\prime \prime}<0$.
\end{conj}

\section{The log-behavior of the sequence $\{a_n(\mu)\}_{n\geq 1}$}

In this section, we study the log-behavior of the sequence
$\{a_n(\mu)\}_{n\geq 1}$. We begin with the log-behavior of the
Bessel zeta functions $\zeta_{\mu}(x)$.

\begin{lemm}\label{besselzetalv} For $\mu >-1$,
 the Bessel zeta function $\zeta_{\mu}(x)$  is log-convex
 for $x>1$.
\end{lemm}
\proof We proceed to show that for $x>1$,
\[
(\log \zeta_{\mu}(x))^{\prime \prime}>0,
\]
or equivalently,
\begin{equation}
  \zeta_{\mu}(x)\zeta_{\mu}^{\prime\prime}(x)-(\zeta_{\mu}^{\prime}(x))^2>0.
\end{equation}
By the convergence of $\zeta_{\mu}(x)$, it is easily seen that
\[
  \zeta_{\mu}^{\prime}(x)=-\sum_{n=1}^{\infty}\frac{\log j_{\mu,n}}{j_{\mu,n}^x}
\]
and
\[
 \zeta_{\mu}^{\prime\prime}(x)=\sum_{n=1}^{\infty}\frac{(\log j_{\mu,n})^2}{j_{\mu,n}^x}.
\]
Hence
\begin{eqnarray*}
  \lefteqn{{\hspace{-0.3cm}\zeta_{\mu}(x)\zeta_{\mu}^{\prime\prime}(x)-(\zeta_{\mu}^{\prime}(x))^2}}\\[6pt]
  &=& \sum_{m=1}^{\infty}\frac{1}{j_{\mu,m}^x}\sum_{n=1}^{\infty}\frac{(\log j_{\mu,n})^2}{j_{\mu,n}^x}
     -\sum_{m=1}^{\infty}\frac{\log j_{\mu,m}}{j_{\mu,m}^x}\sum_{n=1}^{\infty}\frac{\log j_{\mu,n}}{j_{\mu,n}^x}\\[6pt]
  &=& \sum_{n>m\geq 1}\frac{(\log j_{\mu,m})^2+(\log j_{\mu,n})^2-2(\log j_{\mu,m})(\log j_{\mu,n})}{j_{\mu,m}^x j_{\mu,n}^x}\\[6pt]
  &=& \sum_{n>m\geq 1}\frac{(\log j_{\mu,m}-\log j_{\mu,n})^2}{j_{\mu,m}^x j_{\mu,n}^x},
\end{eqnarray*}
which is positive. This completes the proof.
\qed

Setting $f(x)=\zeta_{\mu}(x)$, $x=2n-2$, $y=2n+2$ and $\lambda=1/2$
in the defining relation \eqref{df} of a log-convex function, we
obtain that for $\mu>-1$,
\begin{equation}
  \zeta_{\mu}(2n-2)\zeta_{\mu}(2n+2)>\zeta_{\mu}(2n)^2.
\end{equation}
This yields that the sequence $\{\zeta_{\mu}(2n)\}_{n\geq 1}$ is
log-convex for $\mu>-1$. On the other hand, it is easily checked that
the sequence $\{2^{2n+1}(n+1)!(\mu +1)_n\}_{n\geq 1}$ is  log-convex for
$\mu>-1$.
It is well-known that for two positive log-convex sequences
$\{u_n\}_{n\geq 1}$ and $\{v_n\}_{n\geq 1}$, the sequence
$\{u_nv_n\}_{n\geq 1}$ is also log-convex.  So we arrive at the
following property.

\begin{theo}\label{lcc}
  The sequence $\{a_n(\mu)\}_{n\geq 1}$ is log-convex for $\mu>-1$.
\end{theo}
For $\mu=0$ and  $\mu=1$, Theorem \ref{lcc} gives affirmative
answers to the two conjectures of Amdeberhan, Moll and Vignat
\cite{amd} on the log-convexity of the sequences $\{a_n\}_{n\geq 1}$
and $\{b_n\}_{n\geq 1}$, where $a_n=a_n(1)$ and
$b_n=\frac{1}{2}a_n(0)$.

Next we consider the monotone property of the sequence
$\{\sqrt[n]{a_n(\mu)}\}_{n\geq 1}$ for $\mu>0$.

\begin{theo} For $\mu>0$, the sequence $\{\sqrt[n]{a_n(\mu)}\}_{n\geq 1}$
 is increasing for $n> 4e(\mu+2)^2$
\end{theo}

 To prove the above theorem, we
introduce the function
\[
  \theta_{\mu}(x)=\left(\frac{2}{\mu !}\Gamma\left(\frac{x}{2}\right)
  \Gamma\left(\frac{x}{2}+\mu+1\right)\zeta_{\mu}(x)\right)^{\frac{1}{x}},
\]
which  has the following monotone property.

\begin{theo}\label{lccc}
 For   $\mu\geq 0$, the function $\log \theta_{\mu}(x)$ is
 strictly increasing for $x>8e(\mu+2)^2$.
\end{theo}

\proof Assume that $\mu \geq 0$. To prove the monotone property
 in the theorem, we aim to show that for $x>8e(\mu+2)^2$,
\begin{equation}\label{t1}
  (\log \theta_{\mu}(x))^{\prime}>0.
\end{equation}

Let \begin{equation}\label{hx}
 h(x)=\frac{2}{\mu
!}\Gamma(x/2)\Gamma(x/2+\mu+1)\zeta_{\mu}(x).
\end{equation}
Recalling the definition of $\theta_{\mu}(x)$ as given by
\eqref{thetamu}, we have
\[  \theta_{\mu}(x)=h(x)^{\frac{1}{x}}\]
 and
\[ \log \theta_{\mu}(x)=\frac{1}{x}\log h(x).\] It follows that
\begin{equation}\label{lth}  (\log \theta_{\mu}(x))^{\prime}
  =\frac{1}{x}\left(\frac{h^{\prime}(x)}{h(x)}-\frac{\log h(x)}{x}\right).
\end{equation}

Since $\zeta_{\mu}(x)$ and $\Gamma(x)$ are continuous and
differentiable on $(1,\infty)$, so is $h(x)$. We shall apply the
mean value theorem to $\log h(x)$ on $[2,x]$, where
$x>8e(\mu+2)^2$ and $\mu>-1$.  To this end, we need to show that
$h(2)<1$ and $h(x)>1$ for $\mu>-1$ and $x>8e(\mu+2)^2$.

 Recalling the definition of $h(x)$ as
given by \eqref{hx}, we get
\[
 h(2)=\frac{2}{\mu !}\Gamma (1)\Gamma
(\mu+2)\zeta_{\mu}(2),
\]
where
\[ \zeta_{\mu}(2)=\frac{1}{4(\mu+1)},\]
$ \Gamma(1)=1$  and $ \Gamma(\mu+2)=(\mu+1)!$.  Then
\begin{equation}\label{2}
  h(2)=\frac{2}{\mu !}\cdot (\mu+1)!\cdot \frac{1}{4(\mu+1)},
\end{equation}
so $h(2)<1$.

It remains to show that $h(x)>1$ for $\mu>-1$ and $x>8e(\mu+2)^2$.
Recall that
\begin{equation}\label{ubz}
j_{\mu,
1}<(\mu+1)^{\frac{1}{2}}\left((\mu+2)^{\frac{1}{2}}+1\right),
\end{equation}
for $\mu>-1$, see Chamber \cite{cha}. It follows that for $\mu>-1$,
\begin{equation}\label{uu}
j_{\mu, 1}<2(\mu+2).
\end{equation}
Therefore, we obtain that for $\mu>-1$,
\begin{equation}\label{lbzt}
\zeta_{\mu}(x)=\sum_{n=1}^{\infty}\frac{1}{j_{\mu, n }^x}
>\frac{1}{j_{\mu, 1}^x}>\frac{1}{2^x(\mu+2)^x}.
\end{equation}
On the other hand, it is  known that for $x\geq 0$,
\begin{equation}\label{lbg}
\Gamma(x)>\sqrt{2\pi x}\left(\frac{x}{e}\right)^x,
\end{equation}
see Alzer \cite{alz}. Combining \eqref{lbzt} and \eqref{lbg}, we
deduce that for $x>2$ and $\mu>-1$,
\[
2\Gamma\left(\frac{x}{2}\right)\zeta_{\mu}(x)>2\sqrt{\pi x}
\left(\frac{x}{8e(\mu+2)^2}\right)^{\frac{x}{2}}.
\]
Consequently, for $\mu>-1$ and $x>8e(\mu+2)^2$,  we obtain that
\begin{equation}\label{p1}
2\Gamma\left(\frac{x}{2}\right)\zeta_{\mu}(x)>2\sqrt{\pi x}>1.
\end{equation}
Clearly, for $x>0$ we have
\begin{equation}\label{others}
\frac{\Gamma(x/2+\mu+1)}{{\mu}!}>1.
\end{equation}
In view of \eqref{p1} and \eqref{others}, we find that for $\mu> -1$
and $x>8e(\mu+2)^2$,
\begin{equation}\label{1}
h(x)=\frac{2}{\mu !}\Gamma \left(\frac{x}{2}\right)\Gamma
\left(\frac{x}{2}+\mu+1\right)\zeta_{\mu}(x)>1,
\end{equation}
as claimed.

Next we proceed to prove that there exists $t$ in $(2,x)$ such that
\begin{equation}\label{t2} \frac{h^{\prime}(t)}{h(t)}>\frac{\log
h(x)}{x}.
\end{equation}
By the mean value theorem applied to $\log h(x)$ on $[2,x]$, there exists
$t\in (2,x)$ such that
\begin{equation}\label{4}
\frac{h^{\prime}(t)}{h(t)}=(\log h(t))^{\prime}= \frac{\log
h(x)-\log h(2)}{x-2}.
\end{equation}
On the other hand, we have  shown that $ h(2)<1$ and $h(x)>1$ for
$\mu>-1$ and $x>8e(\mu+2)^2$. Consequently, we have $\log h(2)<0$
and $\log h(x)>0$. Note that for $\mu>-1$ and $x>8e(\mu+2)^2$, we
have $x>2$. Hence
\begin{equation}\label{3}
\frac{\log h(x)}{x}<\frac{\log h(x)-\log h(2)}{x-2}.
\end{equation}
Combining \eqref{4} and \eqref{3}, we obtain \eqref{t2}.

Moreover, it can be shown that
\begin{equation}\label{t3}
  \frac{h^{\prime}(x)}{h(x)}> \frac{h^{\prime}(t)}{h(t)}.
\end{equation}
We claim that for $y> 2$,
\begin{equation}\label{h}
  \left( \frac{h^{\prime}(y)}{h(y)}\right)^{\prime}>0.
\end{equation}
By the definition of $h(x)$ as given by \eqref{hx}, we have
\begin{eqnarray*}
\left( \frac{h^{\prime}(y)}{h(y)}\right)^{\prime}
&=&(\log h(y))^{\prime\prime}\\[5pt]
&=&(\log \Gamma(y/2))^{\prime\prime}
+(\log \Gamma(y/2+\mu+1))^{\prime\prime}+(\log \zeta_{\mu}(x))^{\prime\prime}.
\end{eqnarray*}

It is known  that $(\log \Gamma(y))^{\prime\prime}>0$ for $y>1$ ,
see \cite[Theorem 1.2.5]{and}. Thus, $(\log
\Gamma(y/2))^{\prime\prime}>0$ and $(\log
\Gamma(y/2+\mu+1))^{\prime\prime}>0$ for $y>2$. Moreover, in the
proof of Lemma \ref{besselzetalv}, we have shown that $(\log
\zeta_{\mu}(y))^{\prime\prime}>0$. This proves \eqref{h}. In other
words, $\frac{h^{\prime}(y)}{h(y)}$ is strictly increasing
 for $y>2$. Thus for $2<t<x$, \eqref{t3} holds.

Combining \eqref{t2} and \eqref{t3}, for $\mu>-1$ and
$x>8e(\mu+2)^2$, we find that
\[
\frac{h^{\prime}(x)}{h(x)}-\frac{\log h(x)}{x}
>\frac{h^{\prime}(x)}{h(x)}-\frac{h^{\prime}(t)}{h(t)}> 0.
\]
Hence   \eqref{t1} follows from \eqref{lth}.
 This completes the proof.
\qed

In view of  relation \eqref{a}, it can be checked that
\begin{equation}
  \sqrt[n]{a_n(\mu)}=4\theta_{\mu}(2n)^2.
\end{equation}
Thus Theorem \ref{lccc} implies that for any $\mu\geq 0$
and $n> 4e(\mu+2)^2$, we have $\sqrt[n]{a_n(\mu)}<\sqrt[n+1]{a_{n+1}(\mu)}$ .

For $\mu=1$, it can be verified that
$\sqrt[n]{a_n(1)}<\sqrt[n+1]{a_{n+1}(1)}$ for
 $2 \leq n\leq 108$. In the meantime, for $\mu=1$,
 Theorem  \ref{lccc} states that $\sqrt[n]{a_n(1)}<\sqrt[n+1]{a_{n+1}(1)}$
 for $n>101$. Thus we have the following assertion.

\begin{theo}
The sequence $\{\sqrt[n]{a_n}\}_{n\geq 2}$ is strictly increasing.
\end{theo}

For $\mu=0$, it can be verified that
$\sqrt[n]{a_n(0)}<\sqrt[n+1]{a_{n+1}(0)}$ for $2 \leq n\leq 48$.
Meanwhile, for $\mu=0$, Theorem \ref{lccc} states that
$\sqrt[n]{a_n(0)}<\sqrt[n+1]{a_{n+1}(0)}$ for $n>45$. So we have
$\sqrt[n]{a_n(0)}<\sqrt[n+1]{a_{n+1}(0)}$ for $n\geq 2$. Since
$b_n=\frac{1}{2}a_n(0)$, we have for $n\geq 2$,
\[
 \sqrt[n]{b_n}=\frac{\sqrt[n]{a_n(0)}}{\sqrt[n]{2}}
 < \frac{\sqrt[n+1]{a_{n+1}(0)}}{\sqrt[n+1]{2}}=\sqrt[n+1]{b_{n+1}}.
\]
Thus we have the following monotone property.

\begin{theo}
  The sequence $\{\sqrt[n]{b_n}\}_{n\geq 2}$ is strictly increasing.
\end{theo}

Note that Wang and Zhu \cite{wang3} independently proved the
log-convexity of $\{a_n\}_{n\geq 1}$ and $\{b_n\}_{n\geq 1}$ and the
increasing properties of $\{\sqrt[n]{a_n}\}_{n\geq 1}$ and
$\{\sqrt[n]{b_n}\}_{n\geq 1}$.

\section{The log-behavior of Bell numbers}

%

 In this section, we consider the log-behavior of
 Bell numbers, which are also denoted by $B_n$.
 Recall that the function $B(x)$ is defined by
 \[B(x)=\frac{1}{e}\sum_{k=0}^{\infty}\frac{k^x}{k!} .\]

\begin{lemm}\label{bl}
  The function $B(x)$ is log-convex for $x>1$.
\end{lemm}

\proof
We proceed to show that
\[
(\log B(x))^{\prime\prime}>0,
\]
that is,
\begin{equation}\label{bellflv}
  B(x)B^{\prime\prime}(x)-(B^{\prime}(x))^2> 0.
\end{equation}
For $x\geq 1$, by the convergence of $B(x)$, we have
\[
B^{\prime}(x)=\frac{1}{e}\sum_{n=0}^{\infty}\frac{n^x\log n}{n!}
\]
and
\[
B^{\prime \prime}(x)=\frac{1}{e}\sum_{n=0}^{\infty}\frac{n^x (\log n)^2}{n!}.
\]
Thus, for $x>1$, we have
\begin{eqnarray*}
  \lefteqn{{\hspace{-0.3cm}B(x)B^{\prime\prime}(x)-(B^{\prime}(x))^2}}\\[8pt]
  &=&\frac{1}{e^2}\sum_{m=0}^{\infty}\frac{m^x}{m!}\sum_{n=0}^{\infty}\frac{n^x(\log n)^2}{n!}-
      \frac{1}{e^2}\sum_{m=0}^{\infty}\frac{m^x\log m}{m!}\sum_{n=0}^{\infty}\frac{n^x\log n}{n!}\\[8pt]
  &=&\frac{1}{e^2}\sum_{n>m\geq 0}\frac{m^xn^x}{m!n!}((\log m)^2+(\log n)^2-2\log m\log n)\\[8pt]
  &=&\frac{1}{e^2}\sum_{n>m\geq 0}\frac{m^xn^x}{m!n!}(\log n-\log m)^2,
\end{eqnarray*}
which is positive. This completes the proof.
\qed


We now turn to the log-behavior of the function $B(x)^{1/x}$.

\begin{theo}\label{bm}
$\log B(x)^{1/x}$ is strictly increasing for $x\geq 1$.
\end{theo}

\proof To prove that $\log B(x)^{1/x}$ is strictly increasing, we
aim to show that
\begin{equation}\label{bll}
(\log B(x)^{1/x})^{\prime}>0.
\end{equation}
Since
\[
  (\log B(x)^{1/x})^{\prime}
  =\frac{1}{x}\left(\frac{B^{\prime}(x)}{B(x)}-\frac{\log B(x)}{x}\right),
\]
\eqref{bll} can be rewritten as
\begin{equation}\label{bell44}
  \frac{B^{\prime}(x)}{B(x)}>\frac{\log B(x)}{x}.
\end{equation}
We claim that there exists $t$ in $(1,x)$ such that
\begin{equation}\label{bellll}
\frac{B^{\prime}(t)}{B(t)}>\frac{\log B(x)}{x}.
\end{equation}
Since $B(1)=1$ and $B(x)>1$ for $x>1$, by the mean value theorem
with respect to $\log B(x)$ on $[1,x]$, there exists $t\in (1,x)$
such that
\begin{equation}\label{b1}
  \frac{B^{\prime}(t)}{B(t)}=\frac{\log B(x)-\log B(1)}{x-1}
  =\frac{\log B(x)}{x-1}.
\end{equation}
Since $x>1$, we have \begin{equation}\label{b2}
 \frac{\log
B(x)}{x-1}>\frac{\log B(x)}{x}.
\end{equation}
Combining \eqref{b1} and \eqref{b2}, we obtain \eqref{bellll}.

Next we show that for $x>t>1$,
\begin{equation}\label{bel}
  \frac{B^{\prime}(x)}{B(x)}>\frac{B^{\prime}(t)}{B(t)}.
\end{equation}
In fact, by Lemma \ref{bl}, we see that for $y \geq 1$,
\[
  \left(\frac{B^{\prime}(y)}{B(y)}\right)^{\prime}=(\log B(y))^{\prime\prime}>0.
\]
 This implies that $\frac{B^{\prime}(y)}{B(y)}$ is strictly
 increasing for  $y>1$. This proves \eqref{bel}.

Combining \eqref{bellll} and \eqref{bel}, we obtain \eqref{bell44}.
This completes the proof.
\qed

Since $B(n)=B_n$ whenever $n$ is a positive integer, Theorem
\ref{bm} implies the following monotone property  conjectured by Sun
\cite{sun}.

\begin{coro}\label{bell1}
  The sequence $\{\sqrt[n]{B_n}\}_{n\geq 1}$ is strictly increasing.
\end{coro}

The above property was independently obtained by Wang and Zhu
\cite{wang3} via a different approach. Furthermore, we pose the
following conjecture which implies the conjecture of Sun \cite{sun}
stating that the sequence $\{\sqrt[n]{B_n}\}_{n\geq 1}$ is
log-concave.

\begin{conj}
The function $B(x)^{1/x}$ is log-concave for $x\geq 1$,
that is,  $(\log B(x)^{1/x})^{\prime \prime}<0$ for $x>1$.
\end{conj}

\section{A connection to H\"{o}lder's inequality}

In this section, we give a derivation of
the monotone property of the function $B(x)^{1/x}$ as given in
Theorem \ref{bm} by applying  H\"{o}lder's inequality in probability theory.
In fact, it can be shown the condition $1<x<y$ in Theorem
 \ref{bm} can be relaxed to $0<x<y$.

 Let $Z$ be the discrete random variable with Possion distribution
as given by
\[
P(Z=k)=\frac{1}{e}\frac{1}{k!}.
\]
From Dobinski's formula, it is easily checked that
  $B(x)=E[Z^x]$. H\"{o}lder's inequality states that for
  real-valued random variables $U$, $V$ and positive numbers
  $p$ and $q$ satisfying $\frac{1}{p}+\frac{1}{q}=1$, we have
\[
E[|UV|]\leq E[|U|^p]^{1/p}E[|V|^q]^{1/q},
\]
and the equality holds if and only if either there exist constants
$\alpha, \beta
> 0$ such that $\alpha|U|^p=\beta|V|^q$ or $E[|U|^p]=0$ or
$E[|V|^q]=0$, see, for example, Sachkov \cite{sac}. For $0 < x < y$,
we set $p = y/x$, and set  $U=Z^x$ and $V=1$. It is not hard to see
that in this case H\"{o}lder inequality is strict. Hence we obtain
that
\[
E[Z^x]^{1/x}<E[Z^y]^{1/y},
\]
which can be restated as follows.

\begin{theo}\label{bell2}
For $0<x<y$, we have $B(x)^{1/x}<B(y)^{1/y}$.
\end{theo}

 \vspace{.3cm}

\noindent{\bf Acknowledgments.} We wish to thank the referee for valuable suggestions. This work was
supported by the 973 Project,
the PCSIRT Project, the Doctoral Program Fund of the Ministry of Education,
and the National Science Foundation of China.


\begin{thebibliography}{99}

\bibitem{alz} H. Alzer, Sharp upper and lower bounds for the gamma function, Proc. Royal Soc. Edinburgh A 139 (2009) 709--718.

\bibitem{amd} T. Amdeberhan, V.H. Moll and C. Vignat, A probabilistic interpretation of a sequence related to Narayana polynomials, arXiv:1202.1203.

\bibitem{and} G.E. Andrews, R. Askey and R. Roy, Special Functions, Cambridge University Press, 1999.


\bibitem{art} E. Artin, The Gamma Function, Holt, Rinehart and Winston, New York, 1964.


\bibitem{bel} E.T. Bell, Exponential Numbers, Amer. Math. Monthly 41 (1934) 411--419.

\bibitem{car} L. Carlitz, A sequence of integers related to the Bessel function, Proc. Amer. Math. Soc. 14 (1963) 1--9.

\bibitem{cha} L. G. Chamber, An upper bound for the first zero of Bessel functions, Math. Comp. 38 (1982) 589--591.

%

%
%

\bibitem{las} M. Lasalle, Two integer sequences related to Catalan numbers, J. Combin. Theory Ser. A 119 (2012) 923--935.

\bibitem{luca} F. Luca and P. St\u{a}nic\u{a}, On some conjectures on the monotonicity of some combinatorial sequences, J. Combin. Number Theory 4 (2012) 1--10.

\bibitem{slo} The OEIS Foundation Inc., The On-Line Encyclopedia of Integer Sequences, http://oeis.org.



\bibitem{rota} S.M. Roman and Gian-Carlo Rota, The umbral calculus, Adv. in Math. 27 (1978) 95--188.

\bibitem{rota2} G.-C. Rota, The number of partitions of a set, Amer. Math. Monthly 71 (1964) 498--504.

\bibitem{sac} V.N. Sachkov, Probabilistic Methods in Combinatorial Analysis, Cambridge
University Press, New York, 1997.


\bibitem{smi} D.E. Smith, Source Book in Mathematics, Vols. 1 and 2, New York: Dover, 1959.




\bibitem{sun} Z. Sun, Conjectures involving arithmetical sequences, Numbers Theory: Arithmetic in Shangri-La (eds., S. Kanemitsu, H. Li and J. Liu), Proc. 6th China-Japan Seminar (Shang-hai, August 15-17, 2011), World Sci., Singapore, 2013, pp. 244-258.
%

\bibitem{wang3} Y. Wang and B. Zhu, Proofs of some conjectures on monotonicity of number-theoretic and combinatorial sequences, arXiv:1303.5595.


\end{thebibliography}
\end{document}